

\input amstex
\documentstyle{amsppt}
\magnification=\magstep1
\pagewidth{6.5truein}
\pageheight{8.8truein}
\NoBlackBoxes

\topmatter
\title    Annular and Boundary Reducing Dehn Fillings \endtitle
\author    Cameron McA. Gordon$^{1}$ and Ying-Qing Wu \endauthor
\leftheadtext{C.~M{\fiverm c}A.~Gordon and Y-Q.~Wu}
\rightheadtext{Annular and boundary reducing Dehn fillings}
\address Department of Mathematics, The University of Texas at Austin,
Austin, TX 78712  \endaddress
\email  gordon\@math.utexas.edu \endemail
\address Department of Mathematics, University of Iowa, Iowa City, IA
52242 \endaddress
\email  wu\@math.uiowa.edu \endemail
\thanks  $^1$ Partially supported by NSF grant \#DMS 9626550. \endthanks 
\endtopmatter
 
\document
\define\proof{\demo{Proof}}
\define\endproof{\qed \enddemo}

\redefine\hat{\widehat}

\redefine\bdd{\partial}
\define\Int{\text{\rm Int}}
\baselineskip 16pt
\input epsf.tex
\TagsOnRight
\define\a{\alpha}
\redefine\b{\beta}

\redefine\L{\Lambda}

\head \S 0.  Introduction
\endhead

Surfaces of non-negative Euler characteristic, i.e., spheres, disks,
tori and annuli, play a special role in the theory of 3-dimensional
manifolds.  For example, it is well known that every (compact,
orientable) 3-manifold can be decomposed into canonical pieces by
cutting it along essential surfaces of this kind \cite{K}, \cite{M},
\cite{Bo}, \cite{JS}, \cite{Jo1}.  Also, if (as in \cite{Wu3}) we call a
3-manifold that contains no essential sphere, disk, torus or annulus
{\it simple\/}, then Thurston has shown \cite{T1} that a 3-manifold
$M$ with non-empty boundary is simple if and only if $M$ with its
boundary tori removed has a hyperbolic structure of finite volume with
totally geodesic boundary.  For closed 3-manifolds $M$, the
Geometrization Conjecture \cite{T1} asserts that $M$ is simple if and
only if $M$ is either hyperbolic or belongs to a certain small class
of Seifert fiber spaces.

Because of their importance, a good deal of attention has been
directed at the question of when surfaces of non-negative Euler
characteristic can be created by Dehn filling.  To describe this, let
$M$ be a simple 3-manifold, with a torus boundary component
$\partial_0 M$.  Let $\alpha$ be the isotopy class of an essential
simple loop (or {\it slope\/}) on $\partial_0M$.  Recall that the
manifold obtained from $M$ by {\it $\alpha$-Dehn filling\/} is
$M(\alpha) = M\cup V_\alpha$, where $V_\alpha$ is a solid torus, glued
to $M$ by a homeomorphism between $\partial_0M$ and $\partial
V_\alpha$ which identifies $\alpha$ with the boundary of a meridian
disk of $V_\alpha$.  We are interested in obtaining restrictions on
when $M(\alpha)$ fails to be simple.  Although clearly little can be
said in general about a single Dehn filling, if one considers pairs of
non-simple fillings $M(\alpha)$, $M(\beta)$ then it turns out that the
{\it distance\/} $\Delta (\alpha,\beta)$ between the two slopes
$\alpha$ and $\beta$ (i.e., their minimal geometric intersection
number) is quite small, and hence a given $M$ can have only a small
number of non-simple fillings.  More precisely, if $M(\alpha)$,
$M(\beta)$ contain essential surfaces $F_\alpha,F_\beta$ of
non-negative Euler characteristic, then for each of the ten possible
pairs of homeomorphism classes of $F_\alpha,F_\beta$ one can obtain
upper bounds on $\Delta (\alpha,\beta)$.  In the present paper we deal
with the case where $F_\alpha$ is an annulus and $F_\beta$ is a disk,
and prove the following theorem.

\proclaim{Theorem 0.1} Let $M$ be a simple 3-manifold such that
$M(\alpha)$ is annular and $M(\beta)$ is boundary reducible.  Then
$\Delta (\alpha,\beta)\le 2$.  
\endproclaim

The assumption that $M$ is a simple manifold can be replaced by the
weaker assumption that it is boundary irreducible and anannular, see
Corollary 5.5.  The bound is sharp: infinitely many examples of simple
3-manifolds $M$ with $M(\alpha)$ annular, $M(\beta)$ a solid torus,
and $\Delta(\alpha, \beta) =2$ are given in \cite{MM2}.  See also
\cite{EW}.

Theorem 0.1 completes the determination of the best possible upper
bounds on $\Delta (\alpha,\beta)$ in all ten cases.  These are shown
in Table~1, where $S$, $D$, $A$ and $T$ indicate that the manifold
$M(\alpha)$ or $M(\beta)$ contains an essential sphere, disk, annulus
or torus, respectively.
References for these bounds are: 
$(S,S)$: \cite{GL3} (see also \cite{BZ}); 
$(S,D)$: \cite{Sch}; 
$(S,A)$: \cite{Wu3}; 
$(S,T)$: \cite{Wu1}, \cite{Oh}; 
$(D,D)$: \cite{Wu2}; 
$(D,T)$: \cite{GL4}; 
$(A,A)$, $(A,T)$ and $(T,T)$: \cite{Go}. 
Examples showing that the bounds are best possible can be found in: 
$(S,S)$: \cite{GLi}; $(D, D)$: [Be] and [Ga]; $(S, A)$,
$(D, A)$ and $(D, T)$: [HM]; $(S, T)$: [BZ2]; $(T, A)$ and $(A, A)$:
[GW]; $(T, T)$: [T2] and [Go].

\bigskip
\leavevmode

\centerline{\epsfbox{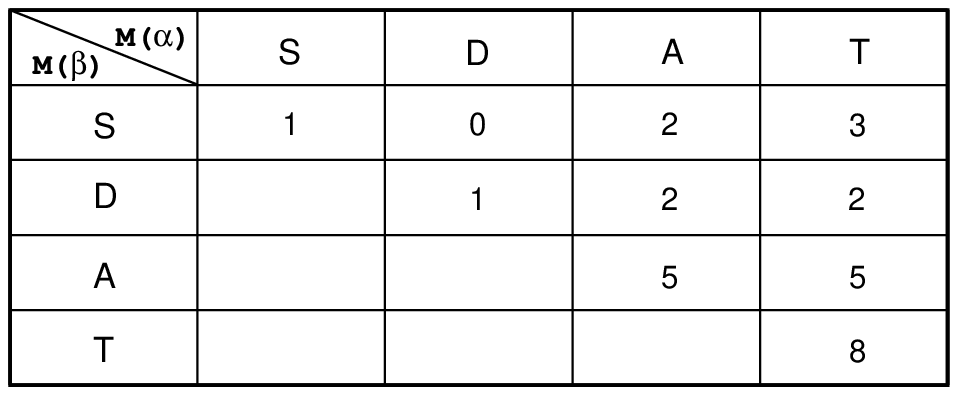}} 
\bigskip
\centerline{Table 1: Upper bounds on $\Delta(\a, \b)$}
\bigskip

Here is a sketch of the proof of Theorem 0.1.  It has been shown by
Qiu [Qiu] that $\Delta \leq 3$, so we assume $\Delta = 3$, and try to
get a contradiction.  Let $A$ and $B$ be an essential annulus and an
essential disk in $M(\a)$ and $M(\b)$, and let $P$ and $Q$ be the
intersection of $A$ and $B$ with $M$, respectively.  Let $p, q$ be the
number of boundary components of $P, Q$ on the torus $\bdd_0M$.
Denote by $K_{\b}$ the core of the Dehn filling solid torus $V_{\b}$
in $M(\b)$.

In Section 2 we consider the special case that $K = K_{\b}$ is a {\it
1-arch\/} knot, which means that it can be isotoped to a union of two
arcs $C_1$ and $C_2$, such that $C_1$ lies on $\bdd M(\b)$, and $C_2$
is disjoint from the compressing disk $B = F_{\b}$ of $\bdd M(\b)$.
In this case the manifold $M(\b)$ is homeomorphic to a manifold $X_C$
obtained by adding a 2-handle to a certain manifold $X$ along a curve
$C$.  This changes a Dehn surgery problem to a handle addition
problem, and we will use a theorem of Eudave-Mu\~noz to show that in
this case the annulus $A = F_{\a}$ can be chosen to intersect the knot
$K_{\a}$ at most twice, that is, $p \leq 2$.

As usual, the intersection of $P \cap Q$ defines graphs $G_A, G_B$ on
$A$ and $B$.  In Section 3 the ``representing all types'' techniques
developed in [GL1--GL4] is modified to suit the case that the
intersection graphs have boundary edges.  It will be proved that when
$\Delta \geq 2$, either $G_A$ represents all types, or $G_B$ contains
a great web.  The first possibility is impossible because it would
lead to a boundary reducing disk of $M(\b)$ which has less
intersection with $K_{\b}$, hence $G_B$ must contain a great web.
This great web is then used in Section 4 to show that if $p \geq 3$
then the knot $K_{\b}$ is a 1-arch knot.  Combined with the result of
Section 2, this proves Theorem 0.1 in the generic case that $p \geq
3$.  Finally in Section 5 the case $p \leq 2$ is ruled out, completing
the proof of Theorem 0.1.

The authors would like to thank John Luecke for helpful conversations.

\head \S 1.  Preliminaries
\endhead

Recall that a 3-manifold $X$ is {\it boundary reducible\/} if its
boundary, denoted by $\bdd X$, is compressible in $X$, in which case a
compressing disk of $\bdd X$ is also called a {\it boundary reducing
disk\/} of $X$.  A surface of non-positive Euler characteristic in $X$
is {\it essential\/} if it is incompressible, $\bdd$-incompressible,
and is not boundary parallel; a sphere (resp.\ disk) is essential if
it is a reducing sphere (resp.\ boundary reducing disk.)

Let $M$ be a simple 3-manifold, with a torus boundary component
$\bdd_0 M$.  Let $\alpha, \beta$ be slopes on $\bdd_0 M$ such that
$M(\alpha)$ is annular and $M(\beta)$ is boundary reducible.  Let $A$
be an essential annulus in $M(\alpha)$, and let $B$ be an essential
disk in $M(\beta)$.  These give rise to a punctured annulus $P = A
\cap M$ and a punctured disk $Q = B \cap M$ in $M$, where $\bdd_0 P =
P \cap \bdd_0 M$ consists of $p$ copies of $\alpha$, and $\bdd_0 Q = Q
\cap \bdd_0 M$ consists of $q$ copies of $\beta$.  We assume that $A,
B, P, Q$ are chosen so that $p$ and $q$ are minimal.  Note that $p, q$
are positive because $M$ is simple.  Now isotope $P$ and $Q$ to
minimize $|P \cap Q|$, the number of components of $P\cap Q$.  Then no
arc component of $P \cap Q$ is boundary parallel in $P$ or $Q$; no
circle component of $P \cap Q$ bounds a disk in $P$ or $Q$; and each
component of $\bdd_0 P$ meets each component of $\bdd_0 Q$ in $\Delta
= \Delta(\a, \b)$ points.

Ruifeng Qiu showed in [Qiu] that if $M$ is a simple manifold, $M(\a)$
is annular and $M(\b)$ is boundary reducible, then $\Delta \leq 3$.
Thus to prove Theorem 0.1, we need only rule out the possibility that
$\Delta = 3$.  {\it In this paper except in Section 3, we will
assume that $\Delta = 3$, and proceed to get a contradiction.}
Results in Section 3 have been proved in a broader setting, so they
can be used in the future.

Regarding the components of $\bdd_0 P$, $\bdd_0 Q$ as fat vertices, we
get graphs $G_A, G_B$ in $A, B$ respectively, where the edges of $G_A$
and $G_B$ are the arc components of $P\cap Q$ that have at least one
endpoint on $\bdd_0 M$.  Let $J = A$ or $B$.  An edge of $G_J$ is an
{\it interior\/} edge if each of its endpoints lies on a vertex of
$G_J$, and a {\it boundary\/} edge if one of its endpoints lies on a
vertex of $G_J$ and the other lies on $\bdd J$.  The faces of $G_J$
correspond in the usual way to components of $J - \Int N(G_J)$.  A
face of $G_J$ is an {\it interior\/} face if it does not meet $\bdd
J$; otherwise it is a {\it boundary\/} face.  Thus the edges in the
boundary of an interior face are interior edges, while the boundary of
a boundary disk face contains some boundary edges.  Denote by $\hat
G_J$ the reduced graph of $G_J$, in which each parallel family of
edges is replaced by a single edge.

Let $u_1, \ldots, u_p$ be the vertices of $G_A$, labeled successively
when traveling along the Dehn filling solid torus $V_{\a}$.  Each
$u_i$ is given a sign according to whether $V_{\a}$ passes $A$ from
the positive side or negative side at this vertex.  Two vertices $u_i,
u_j$ are {\it parallel\/} if they have the same sign, otherwise they
are {\it antiparallel}.  The vertices $v_1, \ldots, v_q$ of $G_B$ are
labeled and signed similarly.

If $e$ is an edge of $G_A$ with an endpoint on $u_i$, then the
endpoint is labeled $j$ if it is on $\bdd u_i \cap \bdd v_j$.  Thus
when going around $\bdd u_i$, the labels of the edge endpoints appear
as $1, 2, ... , q$ repeated $\Delta$ times.  The edge endpoints of
$G_B$ are labeled similarly.

A cycle in $G_A$ or $G_B$ is a {\it Scharlemann cycle\/} if it bounds
a disk with interior disjoint from the graph, and all the edges in the
cycle have the same pair of labels $\{i, i+1\}$ at their two
endpoints, called the {\it label pair\/} of the Scharlemann cycle.  A
pair of edges $\{e_1, e_2\}$ is an {\it extended Scharlemann cycle\/}
if there is a Scharlemann cycle $\{e'_1, e'_2\}$ such that $e_i$ is
parallel and adjacent to $e'_i$. 

We use $N(X)$ to denote a regular neighborhood of a subset $X$ in a
given manifold.  

\proclaim{Lemma 1.1}  {\rm (Properties of $G_A$.)}

(1) {\rm (The Parity Rule)} An edge connects parallel vertices on
$G_A$ if and only if it connects antiparallel vertices on $G_B$.

(2) $G_A$ does not have $q$ parallel interior edges.

(3) $G_A$ contains no Scharlemann cycles.

(4) Each label $x \in \{1, ..., q\}$ appears at most once among the
endpoints of a family $\Cal E$ of parallel edges in $G_A$ connecting
parallel vertices; in particular, $\Cal E$ contains at most $q/2$
edges.

(5) No pair of edges are parallel on both $G_A$ and $G_B$.
\endproclaim

\proof
(1) This is on [CGLS, Page 279].  

(2) If $G_A$ contains $q$ parallel interior edges, then the core of
the Dehn filling solid torus in $M(\b)$ would be a cable knot, in
which case $M$ contains an essential annulus, contradicting the
assumption.  See the proof of [GLi, Proposition 1.3].

(3) This follows from [CGLS, Lemma 2.5.2].  

(4) If some label appears twice among the endpoints of a family of
parallel edges connecting a pair of parallel vertices, then there is a
Scharlemann cycle among this family, contradicting (3).  See [CGLS,
Lemma 2.6.6].

(5) If a pair of edges are parallel on both $G_A$ and $G_B$, then they
cut off a disk on each of $P$ and $Q$, whose union is an annulus in
$M$, which is essential because its intersection with $\bdd_0 M$ is a
curve intersecting $\a$ at a single point.  This contradicts the
assumption that $M$ is simple.  
\endproof

\proclaim{Lemma 1.2} {\rm (Properties of $G_B$.)}

(1) If $G_B$ has a Scharlemann cycle, then $A$ is a separating
annulus, and $p$ is even.  Moreover, the subgraph of $G_A$ consisting
of the edges of the Scharlemann cycle and their vertices is not
contained in a subdisk of $A$.  

(2) If $p>2$, then $G_B$ has no extended Scharlemann cycle.  Any two
Scharlemann cycles of $G_B$ have the same label pair.  
\endproclaim

\proof (1) This follows from the proof of [CGLS, Lemma 2.5.2].  It was
shown that using the disk bounded by the Scharlemann cycle one can
find another annulus $A'$ in $M(\a)$ which has fewer intersections
with the Dehn filling solid torus, and is cobordant to $A$, so if $A$
were nonseparating then $A'$ would still be essential, which would
contradict the minimality of $p$.  If the subgraph $G$ consisting of
the edges of a Scharlemann cycle and their end vertices is contained
in a disk in $A$ then $A'\cup A$ bounds a connected sum of $A\times I$
and a lens space, so $A$ being essential implies that $A'$ is
essential, which again contradicts the minimality of $p$.

(2) This is [Wu3, Lemma 5.4(2) -- (3)].  If $G_B$ has an extended
Scharlemann cycle or two Scharlemann cycles with distinct label pairs,
then one can find another essential annulus in $M(\a)$ having fewer
intersection with $K_{\a}$, which would contradict the minimality of
$p$.  \endproof

\head \S 2.  1-arch knots
\endhead

Let $K = K_{\b}$ be the core of the Dehn filling solid torus in
$M(\b)$.  The knot $K$ is a {\it 1-arch\/} knot (with respect to
$B$) if $K$ is isotopic to a union of two arcs $C_1$ and $C_2$, such
that $C_1$ lies on $\bdd M(\b)$, and $C_2$ is disjoint from a
compressing disk $B$ of $\bdd M(\b)$.  

Fix an orientation of $K$ so that when traveling along $K$ with this
orientation one meets the fat vertices $v_1, \ldots, v_q$
successively.  Let $K[i]$ be the point $K \cap v_i$, and for $i\neq
j$, let $K[i, j]$ be the oriented arc segment of $K$ starting from
$K[i]$ and ending at $K[j]$.  Thus $K = K[i, j] \cup K[j, i]$.

\proclaim{Lemma 2.1} If $G_A$ contains $q$ parallel boundary edges,
then $K = K_{\b}$ is a 1-arch knot.  \endproclaim

\proof Without loss of generality we may assume that the interior
endpoints of the parallel boundary edges $e_1, \ldots, e_q$ are
successively labeled $1, \ldots, q$.  Let $D$ be the disk on $P$ cut
off by $e_1$ and $e_q$.  Then $D$ can be extended into the Dehn
filling solid torus $N(K)$ to get a disk $D'$ in $M(\b)$ such that
$\bdd D' = e'_1 \cup K[1,q] \cup e'_q \cup C_1$, where $e'_i$ is an
arc on $B$ containing $e_i$, connecting $K[i]$ to the endpoint of
$e_i$ on $\bdd A$, and $C_1 = D \cap \bdd A$ lies on $\bdd M(\b)$.
Now $K$ is isotopic to $C_1 \cup (e'_1 \cup K[q, 1] \cup e'_q)$ via
the disk $D'$.  Let $C_2 = e'_1 \cup K[q, 1] \cup e'_q$.  After a
slight isotopy one can make $B$ disjoint from $C_2$, as desired.
\endproof

\proclaim{Lemma 2.2} Suppose $\hat G_A$ has a vertex $u$ of valency 4,
such that one of the four edges of $\hat G_A$ incident to $u$ is a
boundary edge, and the two edges adjacent to it are interior edges.
Then either $G_B$ contains a Scharlemann cycle, or $K = K_{\b}$ is a
1-arch knot.  \endproclaim

\proof Let $\hat e_1, \hat e_2, \hat e_3, \hat e_4$ be the four edges
of $\hat G_A$ incident to $u$, and assume that $\hat e_2$ is a
boundary edge.  By Lemmas 2.1 and 1.1(2) we may assume that each $\hat
e_i$ represents at most $q-1$ parallel edges of $G_A$.  Now each label
appears at most twice among the endpoints at $u$ of edges represented
by $\hat e_2$ or $\hat e_4$, hence all labels appear on endpoints at
$u$ of edges represented by $\hat e_1$ or $\hat e_3$.  Suppose $\hat
e_1, \hat e_3$ connect $u$ to $u'$ and $u''$, respectively.  If both
$u'$ and $u''$ are antiparallel to $u$, then by the parity rule each
vertex $v$ on $G_B$ is incident to an edge connecting it to a parallel
vertex, with label $u$ at its endpoint at $v$.  By [CGLS, Lemmas 2.6.3
and 2.6.2] this implies that $G_B$ contains a great $u$-cycle, hence a
Scharlemann cycle, and we are done.  Also, notice that $u'$ and $u''$
cannot both be parallel to $u$, otherwise by Lemma 1.1(4) each of
$\hat e_1$ and $\hat e_3$ represents at most $q/2$ edges, and since
each of $\hat e_2$ and $\hat e_4$ represents at most $q-1$ edges, this
would contradict the fact that the total valency of $u$ in $G_A$ is
$\Delta q = 3 q$.  (This also takes care of the case that $\hat e_1 =
\hat e_3$ is a loop at $u$.)  Therefore, we may assume that $u'$ is
parallel to $u$, and $u''$ is antiparallel to $u$.  Since the total
number of edges represented by $\hat e_1 \cup \hat e_2 \cup \hat e_3$
is more than $2q$, we can choose $2 q$ successive edges at $u$,
forming a subgraph as shown in Figure 2.1.  One can now use [Wu2,
Lemma 2.2] and the proof of [Wu2, Lemma 3.4] to show that there is a
disk $D$ in $M(\b)$ with $\bdd D = K[1,q] \cup \a_1 \cup C_1 \cup
\a_2$, where $\a_1$ and $\a_2$ are arcs on the compressing disk $B$
connecting $K[1]$ and $K[q]$ to $\bdd B$, and $C_1$ is an arc on $\bdd
M(\b)$.  As in the proof of Lemma 2.1, this implies that $K$ is a
1-arch knot.  \endproof

\bigskip
\leavevmode

\centerline{\epsfbox{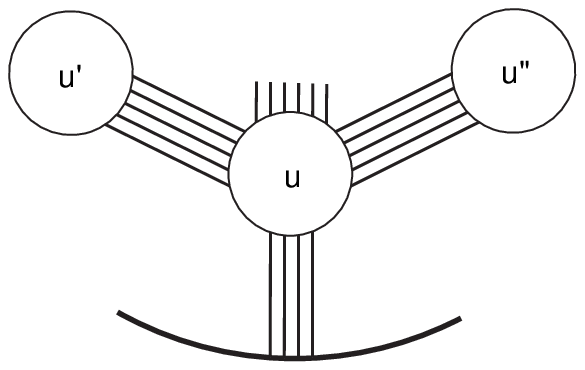}} 
\bigskip
\centerline{Figure 2.1}
\bigskip

We need the following result of Eudave-Mun\~oz in the proof of
Proposition 2.4.  If $C$ is a simple loop on the boundary of a
3-manifold $X$, denote by $X_C$ the manifold obtained by adding a
2-handle to $X$ along the curve $C$.

\proclaim{Lemma 2.3}
Let $X$ be an irreducible, orientable 3-manifold with $\bdd X$
compressible, and $C$ a simple closed curve on $\bdd X$ such that
$\bdd X - C$ is incompressible.  Suppose $X_C$ contains an essential
annulus $A'$.  Then it contains an essential annulus $A$ which
intersects the attached 2-handle in at most two disks.  Furthermore,
if $A'$ is nonseparating, then $A$ can be chosen to be disjoint from
the attached 2-handle.
\endproclaim

\proof This is essentially [Eu, Theorem 1].  The theorem there says
that under the above assumption, either one can find $A$ to be
disjoint from the attached 2-handle, or after sliding the cocore
$\sigma$ of the attached 2-handle over itself to get a 1-complex
$\tau$, one can find an essential annulus $A$ which intersects $\tau$
at a single point.  Moreover, if $A'$ is nonseparating, then $A$ is
disjoint from $\tau$ (see also [Jo2, Sch]).  Sliding $\tau$ back to
$\sigma$, we see that $A$ is isotopic to an annulus intersecting
$\sigma$ at most twice.  See also the remarks after the
statement of Theorem 2 in [Eu].
\endproof

\proclaim{Proposition 2.4} If $K = K_{\b}$ is a 1-arch knot in
$M(\b)$, then $p\leq 2$, and $A$ is a separating annulus in $M(\a)$.
\endproclaim

\proof The first part of the proof here is the same as that in the
proof of [Wu2, Proposition 1].  Suppose $K$ is isotopic to $C_1 \cup
C_2$ as in the definition of 1-arch knot.  Let $Y$ be the manifold
obtained by adding a 1-handle $H_1$ to $M(\b)$ along two disks
centered at $\bdd C_1$, and let $C$ be a simple closed curve on $\bdd
Y$ obtained by taking the union of $C_1 \cap \bdd Y$ and an arc on
$\bdd H_1$.  Let $K'$ be the union of $C_2$ and the core of the
1-handle $H_1$.  Then after adding a 2-handle $H_2$ to $Y$ along $C$
the 1-handle and the 2-handle cancel each other and we get a manifold
$M'$ homeomorphic to the original manifold $M(\b)$, with the knot $K$
identified to $K'$; hence we have a homeomorphism of pairs $(M(\b), K)
\cong (M', K')$.  Let $W$ (denoted by $Q$ in [Wu2]) be the manifold
obtained from $Y$ by Dehn surgery on $K'$ along the slope $\a$.  Then
$M(\a)$ is homeomorphic to the manifold $W_C$ obtained by adding the
2-handle $H_2$ to $W$ along the curve $C$.  It was shown in [Wu2,
Lemmas 1.2 and 1.3] that $\bdd W$ is compressible, and $\bdd W - C$ is
incompressible in $W$ when $\Delta \geq 2$.

Since $M(\b)$ is $\bdd$-reducible, by [Sch] the manifold $M(\a) = W_C$
is irreducible.  This implies that $W$ is irreducible because a
reducing sphere in a manifold always remains a reducing sphere after
2-handle additions.  Therefore we can apply Lemma 2.3 and conclude
that there is an essential annulus $A$ in $M(\a) = W_C$ intersecting
the attached 2-handle $H_2$ in $n\leq 2$ disks; moreover, if $A$ is
nonseparating, then it is disjoint from $H_2$.  Our goal is to show
that $A$ also intersects the knot $K_{\a}$ in $M(\a)$ in two or zero
points, respectively.  

We assume $n=2$, the cases $n=0$ or $1$ are similar.  Let $D_1, D_2$
be the disks $A \cap H_2$, and let $F$ be the twice punctured annulus
$A - \Int (D_1 \cup D_2)$ in $W$.  A meridian disk $D$ of the 1-handle
$H_1$ gives rise to a nonseparating essential annulus $D_0 = D \cap X$
in the manifold $X = Y - \Int N(K') = W - \Int N(K_{\a})$.  Let $F_0 =
F \cap X$.  Form intersection graphs $G_D$ and $G_F$ in the usual way,
i.e, $G_F$ has $F\cap N(K_{\a})$ as fat vertices, $G_D$ has a single
vertex $D\cap N(K_{\b})$, and the edges of $G_D$ and $G_F$ are the arc
components of $D_0 \cap F_0$ which has at least one endpoint on the
fat vertices.  See Figure 2.2.  Choose $A$ and $F$ so that $A_0$
intersects $F_0$ minimally.  Then each fat vertex of $G_F$ has valency
$\Delta = 3$, and the only vertex $x$ of $G_D$ has valency $3t$, where
$t$ is the number of vertices of $G_F$.  Note that, since $A$ is
essential in $M(\a)$, we have $t \geq p$.  As usual, there are no
trivial loops.  Hence each edge of $G_D$ connects $x$ to $\bdd D$.

Each of $\bdd D_1$ and $\bdd D_2$ intersects $\bdd D$ at a single
point, which we denote by $z_1, z_2$, as indicated by the dark dots in
Figure 2.2(a) and (b).  They divide $\bdd D$ into two arcs $\a_1,
\a_2$, one of which, say $\a_1$, lies in $N(C)$, which is the
attaching region of the two handle $H_2$ above.  Hence the interior of
$\a_1$ is disjoint from $\bdd F$.  It follows that all the endpoints
of the edges of $G_D$ on $\bdd D$ lie on the arc $\a_2$, as shown in
Figure 2.2(a).

\bigskip
\leavevmode

\centerline{\epsfbox{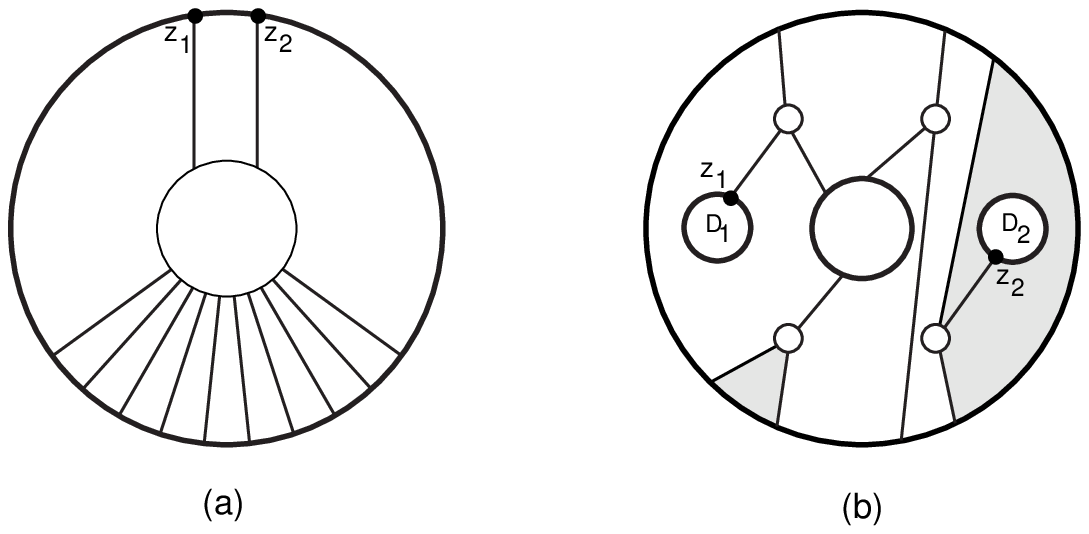}} 
\bigskip
\centerline{Figure 2.2}
\bigskip

Now suppose $G_F$ has $t\geq 3$ vertices.  Since each of $\bdd D_1$
and $\bdd D_2$ is adjacent to at most one edge, there is a vertex of
$G_F$ with two edges connecting it to the same component of $\bdd A
\subset \bdd F$.  An outermost such vertex has a pair of edges $a_1,
a_2$ on $G_F$, cutting off a region $B_1$ on $F_0$ which is either a
disk, or a once punctured disk containing one of $\bdd D_1, \bdd D_2$,
as shown by the two shaded regions in Figure 2.2(b).  (If $B_1$
contains both $\bdd D_i$, choose another outermost vertex.)  They also
cut off a disk $B_2$ on $D_0$ which, by the property in the last
paragraph, has boundary disjoint from $\a_1$.  Therefore, $B_1\cup
B_2$ is either an annulus or a once punctured annulus in $X$, with one
boundary component $\gamma$ a curve on $\bdd N(K)$, another a curve on
$\bdd X - \bdd N(K)$ disjoint from $C$, and a possible third curve
parallel to $C$.  After capping off the last component by a disk in
the attached 2-handle $H_2$, the surface becomes an annulus in $X_C$.
However, since $X_C = M$, and since the boundary component $\gamma$ of
the annulus on the torus $\bdd N(K) = \bdd_0 M$ is an essential curve,
(essential because $\gamma$ is the union of an arc in $\alpha$ and an
arc in $\beta$, and $\alpha$ intersects $\beta$ minimally,) this
contradicts the fact that $M$ is $\bdd$-irreducible and anannular.

When $n=0$, since $t \geq p > 0$, there is a pair of edges which are
parallel on both graphs $G_D$ and$G_F$.  As shown above, this would
give rise to an essential annulus in $M$, which would contradict the
simplicity of $M$.  Hence this case does not happen.  In particular,
this and Lemma 2.3 show that $M(\a)$ cannot contain a nonseparating
annulus.  \endproof

\head \S 3. Representing types
\endhead

Denote by ${\bold q} = \{1, ..., q\}$ the set of labels of the vertices
of $G_B$.  We have the concept of a {\it {\bf q}-type\/} etc.\ from
[GL1].  An interior face of $G_A$ {\it represents\/} a {\bf q}-type
$\tau$ if it is a disk and represents $\tau$ in the sense of [GL1].
We say $G_A$ {\it represents \/} $\tau$ if some interior face of $G_A$
represents $\tau$.

\proclaim{Theorem 3.1} $G_A$ does not represent all {\bf q}-types. 
\endproclaim

\proof See [GL4, Proof of Theorem 2.2].  The proof works for any
essential surface $F$ in $M(\a)$ (in [GL4] $F$ was a torus).  A set of
representatives of all {\bf q}-types contains a set $\Cal D$ of
interior faces of $G_F$ which can be used to surger $Q$ tubed along
the annuli corresponding to the corners of the faces in $\Cal D$,
contradicting the minimality of $q$.  \endproof

A {\it web\/} in $G_B$ is a non-empty connected subgraph $\Lambda$ of
$G_B$ such that all the vertices of $\Lambda$ have the same sign, and
such that there are at most $p$ edge endpoints at vertices of
$\Lambda$ which are not endpoints of edges in $\Lambda$.  Note that a
web may have boundary edges.  

Let $U$ be a component of $B - N(\Lambda)$ that meets $\bdd B$.
Then $D = B - U$ is a {\it disk bounded by $\Lambda$}.  $\Lambda$ is a
{\it great web\/} if there is a disk bounded by $\Lambda$ such that
$\Lambda$ contains all the edges of $G_B$ that lie in $D$.  

\remark{Remark} If there are no boundary edges, then these definitions
coincide with those in [GL2, Section 2].  The following is the analog
in our present setting of [GL2, Theorem 2.3].
\endremark

\proclaim{Theorem 3.2}  Suppose $\Delta \geq 2$.  Let $L$ be a subset
of {\bf q}, and $\tau$ be a non-trivial $L$-type such that

(i)  all elements of $C(\tau)$ have the same sign, and 

(ii) all elements of $A(\tau)$ have the same sign.

If $G_A(L)$ does not represent $\tau$ then $G_B$ contains a web
$\Lambda$ such that the set of vertices of $\Lambda$ is a subset of
either $C(\tau)$ or $A(\tau)$.  
\endproclaim

\proof
Regard $G_A(L)$ as a graph in $S^2$, by capping off the boundary
components of $A$ with two additional fat vertices $v_1, v_2$.  

Define a directed graph $\Gamma = \Gamma(\tau)$ as follows.  The
vertices of $\Gamma$ are the fat vertices of $G_A(L)$ plus
$v_1, v_2$, together with dual vertices of $G_A(L)$ (one in the
interior of each face of $G_A(L)$.)  The edges of $\Gamma$ join each
dual vertex to the fat vertices in the boundary of the corresponding
face.  The edges of $\Gamma$ are oriented as follows:  If an edge $e$
has an endpoint on a vertex of $G_A(L)$, then it is oriented according
to the type $\tau$, (as in [GL1], where $\Gamma$ is denoted by
$\Gamma(T)^*$); if $e$ has an endpoint on $v_1$ or $v_2$, orient $e$
so that no dual vertex in a boundary face of $G_A(L)$ is a sink or
source of $\Gamma$.  See Figure 3.1.

\bigskip
\leavevmode

\centerline{\epsfbox{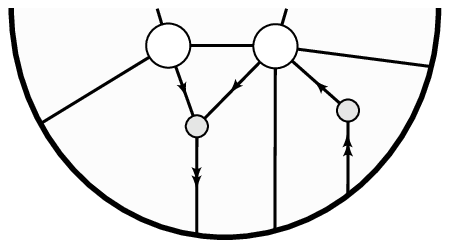}} 
\bigskip
\centerline{Figure 3.1}
\bigskip

By Glass' index formula (see [Gl]) applied to $\Gamma$, we have
$$ \sum_{\text{vertices}} I(v) + \sum _{\text{faces}} I(f) = \chi(S^2)
= 2.$$ 
Assume $G_A(L)$ does not represent $\tau$.  Then no dual vertex in
$\Gamma$ is a sink or source.  Hence
$$\sum _{v \text{ dual}} I(v) \leq 0.$$
Let 
$$ 
\align
&c(\tau) = \text{\# clockwise switches (= \# anticlockwise switches)
of $\tau$} \\
&c(v_i)  = \text{\# clockwise switches (= \# anticlockwise switches)
at $v_i$, $i=1,2$}
\endalign
$$

For $v$ a vertex of $G_A(L)$, we have 
$$ I(v) = 1- \Delta c(\tau).$$
Also, 
$$ I(v_i) = 1- c(v_i), \qquad i = 1,2.$$
Therefore, the number of switch edges, including all switch boundary
edges, is at least 
$$ 
\align
\sum I(f) &\geq 2 + p( \Delta c(\tau) - 1) + (c(v_1) - 1) + (c(v_2)
-1) \\
&= p( \Delta c(\tau) - 1) + c(v_1)+ c(v_2) \tag *
\endalign
$$
Since the number of switch edge endpoints is twice the number of
switch edges, this is also a lower bound for the number of (say)
clockwise switch edge endpoints.  The total number of clockwise
switches is $p \Delta c(\tau) + c(v_1) + c(v_2)$, so the number of
clockwise switches that are not endpoints of clockwise switch edges is
at most $p$.

Since $\Delta \geq 2$, the right hand side of (*) is positive.  Let
$\Lambda$ be a component of the subgraph of $G_B$ consisting of the
edges corresponding to the clockwise switch edges of $G_A(L)$.  Then
at most $p$ edge endpoints at vertices of $\Lambda$ do not belong to
edges of $\Lambda$.  Thus $\Lambda$ is a web, as described.
\endproof

The following is the analog of [GL2, Theorem 2.5].

\proclaim{Theorem 3.3}  Suppose $\Delta \geq 2$.  Let $\Lambda$ be
either (i) a web in $G_B$, or (ii) the empty set.  In case (i), let
$D$ be a disk bounded by $\Lambda$, and in case (ii), let $D = B$.
Let $L$ be the set of vertices of $G_B - \Lambda$ that lie in $D$.
Then either $G_B$ contains a great web or $G_A(L)$ represents all
$L$-types.  
\endproclaim

\proof
Basically, this follows from the proof of [GL2, Theorem 2.5].  We
indicate briefly how this goes.  

We prove the result by induction on $|L|$.

Let $\tau$ be an $L$-type.  We show that if $G_A(L)$ does not
represent $\tau$ then $G_B$ contains a great web.  There are two
cases.

CASE 1.  {\it $\tau$ is trivial.}  Proceed as in [GL2, Proof of Theorem
2.5].  Let $\hat {\Lambda}$ be a component of the subgraph of $G_B$
consisting of vertices $J$, all interior edges with both endpoints on
vertices in $J$, and all boundary edges with one endpoint on a vertex
in $J$.  

(a)  $\hat{\Lambda}$ is a web.  Argue as in [GL2], with ``faces''
meaning ``interior faces''.  

(b)  $\hat{\Lambda}$ is not a web.  Again the argument in [GL2]
remains valid.  More precisely, since $\hat{\Lambda}$ is not a web,
there are more than $p$ edges of $G_B$ connecting a vertex of
$\hat{\Lambda}$ to an antiparallel vertex.  Let $\Sigma$ be the
subgraph of $G_A$ consisting of the vertices of $G_A$ together with
those edges.  Note that these are interior edges of $G_A$, connecting
parallel vertices.  Applying Euler's formula to $\Omega$, a graph in
$A$, gives
$$ V-E+\sum \chi(f) = 0.$$
Therefore
$$\sum\chi(f) = E - V > p - p = 0.$$
Hence $\Sigma$ has a disk face, which must be an interior face.  This
face then contains a face of $G_A(L)$ representing the trivial type.

CASE 2.  {\it $\tau$ is non-trivial.}  Here the argument in [GL2] goes
through essentially without change, (using Theorem 3.2), where we
always interpret ``face'' as ``interior face''.  In particular, [GL2,
Lemmas 2.4 and 2.6] carry over in this way.
\endproof

Theorems 3.1 and 3.3 (in case (ii)) imply:

\proclaim{Corollary 3.4} If $\Delta \geq 2$, then $G_B$ contains a
great web.  \qed \endproclaim

\head \S 4.  The generic case
\endhead

Let $\Lambda$ be a great web in $G_B$ given by Corollary 3.4, and let
$D$ be a disk bounded by $\Lambda$ with the property in the definition
of a great web.  Let $x$ be a label of the vertices of $G_A$, and let
$\Lambda_x$ be the subgraph of $\Lambda$ consisting of all vertices of
$\Lambda$ and all edges in $\Lambda$ with an endpoint labeled $x$.
Let $V$ be the number of vertices of $\Lambda$.  A {\it ghost endpoint
of $\Lambda$\/} is an endpoint, at a vertex of $\Lambda$, of an edge
of $G_B$ which does not belong to $\Lambda$.  A {\it ghost endpoint of
$\Lambda_x$\/} is a ghost endpoint of $\Lambda$ labeled $x$.  (It is
called a ghost label in [GL2].)  By the definition of a web, $\Lambda$
has at most $p$ ghost endpoints.

By a {\it monogon\/} we mean a disk face with one edge in its
boundary, and by a {\it bigon\/} we mean a disk face with two edges in
its boundary.

\proclaim{Lemma 4.1} {\rm (Cf.\ [GL2, Lemma 4.2])}  If $\Lambda_x$ has
at least $3V -2$ edges then $\Lambda_x$ contains a bigon in $D$.
\endproclaim

\proof
Let $\Omega$ be the graph in $S^2$ obtained from $\Lambda_x$ by
regarding $\bdd B$ as a vertex.  Then $\Omega$ has $V+1$ vertices, $E$
edges (= number of edges of $\Lambda_x$), and the faces of $\Omega$
are the faces of $\Lambda_x$ in $D$ together with an additional face
$f_0$.  Note that $f_0$ is not a monogon.  Suppose $\Lambda_x$
contains no bigon in $D$.  

First suppose $f_0$ is not a bigon.  Then $2E \geq 3F$, where $F =
\sum \chi(f)$ summed over all faces of $\Omega$.  Also, 
$$ (V+1) - E + F = 2.$$
Hence 
$$1 = V - E + F \leq V - E + \frac{2E}3 = V - \frac E3,$$
giving $3 \leq 3V - E$, i.e, $E \leq 3V -3$, contrary to 
assumption.  

Now suppose $f_0$ is a bigon; see Figure 4.1.  Then $\Lambda_x$ has at
most one ghost endpoint.  Therefore $E \geq 3V -1$.  Also, Since
$\Lambda_x$ has no bigon in $D$, we have
$$2E \geq 3F - 1.$$  
Hence, as before, 
$$ 1 = V - E + F \leq V - E + \frac {2E + 1} 3 = V - \frac E3 + \frac
13.$$
Therefore $3 \leq 3V - E + 1$, implying $E \leq 3 V - 2$, a contradiction.
\endproof

\bigskip
\leavevmode

\centerline{\epsfbox{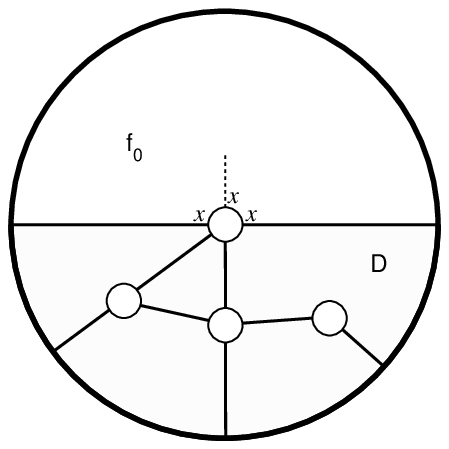}}
\bigskip
\centerline{Figure 4.1}
\bigskip

\remark{Remark} One can show that the conclusion of Lemma 4.1 still
holds if we only assume that $\Lambda_x$ has at least $3V-3$ edges,
but Lemma 4.1 will suffice for our purposes.  
\endremark

\proclaim{Lemma 4.2}  $\Lambda_x$ contains a bigon in $D$ for at least
$2p/3$ labels $x$.  
\endproclaim

\proof (Cf. [GL2, Theorem 4.3]).  By Lemma 4.1, if $\Lambda_x$ does
not contain a bigon in $D$ then $\Lambda_x$ has at most $3V-3$ edges.
Since the vertices of $\Lambda_x$ are all parallel, by the parity rule
no edge of $\Lambda_x$ has both endpoints labeled $x$, so among the
endpoints of edges of $\Lambda_x$, at most $3V-3$ are labeled $x$.
Since $\Delta = 3$, this means that $\Lambda_x$ has at least 3 ghost
endpoints.  Since the total number of ghost endpoints in $\Lambda$ is
at most $p$ by the definition of a great web, there can be at most
$p/3$ such labels $x$.  Hence for at least $2p/3$ labels $x$,
$\Lambda_x$ does contain a bigon in $D$.  \endproof

Note that a boundary bigon in $\Lambda_x$ gives $p+1$ parallel
boundary edges in $G_B$.  

In the remainder of this section, {\it we assume that $p\geq 3$.}

\proclaim{Lemma 4.3}
For some label $x$, $\L_x$ contains a boundary bigon.
\endproclaim

\proof A bigon face of $\Lambda_x$ in $D$ is either a boundary bigon
or an interior bigon.  The latter is either an order 2 Scharlemann
cycle in $G_B$, or contains an extended Scharlemann cycle.  The second
is impossible by Lemma 1.2(2).  When $p$ is odd, the first is also
impossible (Lemma 1.2(1)), and when $p$ is even, any two Scharlemann
cycles have the same label pair (Lemma 1.2(2)).  Hence, by Lemma 4.2,
the number of labels $x$ such that $\Lambda_x$ contains a boundary
bigon is at least
\medskip
$$
\cases
2p/3 \geq 2 \times 3 / 3 = 2, &\text{$p$ odd;} \\
2p/3-2 \geq 2\times 4/3 - 2 = 2/3, &\text{$p$ even.}
\endcases
$$
\medskip
\noindent
Hence there is at least one label $x$ with the stated property.
\endproof

\proclaim{Corollary 4.4} (a) Every vertex of $G_A$ has a boundary edge
incident to it.

(b) $G_A$ has a vertex with two non-parallel boundary edges.
\endproclaim

\proof (a) By Lemma 4.3, $\L_x$ contains a boundary bigon for some
$x$, which gives rise to $p+1$ parallel boundary edges in $G_B$.
Hence each label of ${\bold p} = \{ 1, ..., p\}$ appears at the endpoint
of some boundary edge of $G_B$, and the result follows.

(b) By Lemma 1.1(5), the two boundary edges in a boundary bigon of
$\L_x$ are nonparallel on $G_A$.  
\endproof

\proclaim{Lemma 4.5} $\hat{G}_A$ has no vertex of valency at
most 3.  \endproclaim

\proof $G_A$ has at most $q-1$ parallel interior edges by Lemma
1.1(2), and at most $q-1$ parallel boundary edges by Lemma 2.1,
Proposition 2.4, and the assumption that $p\geq 3$.  Since the total
valency of each vertex of $G_A$ is $3q$, the result follows.
\endproof

\proclaim{Corollary 4.6} $\hat G_A$ has no vertex with two boundary
edges going to the same component of $\bdd A$.
\endproclaim

\proof
Consider an outermost such vertex, with $E$ the corresponding subdisk
of $A$.  Doubling $E$ along the two boundary edges in question and
applying [CGLS, Lemma 2.6.5] gives a vertex in the interior of $E$ of
valency at most 3, contradicting Lemma 4.5.
\endproof

\proclaim{Lemma 4.7} $\hat G_A$ has a vertex $v$ of valency 4, such
that no two boundary edges of $\hat G_A$ at $v$ have endpoints
adjacent on $\bdd v$.  \endproclaim

\proof By Corollaries 4.4(b) and 4.6, $\hat G_A$ has at least one
vertex with two boundary edges going to different components of $\bdd
A$.  Cut $A$ along all such pairs of edges; we get a certain number
($\geq 1$) of disk regions.  If there are no vertices in the interior
of any of these regions, then every vertex of $\hat G_A$ satisfies the
conclusion of the lemma; (recall that there is no vertex of valency
$\leq 3$ by Lemma 4.5).  So consider a region with a non-zero number
of vertices in its interior; see Figure 4.2(a).  Note that each vertex
$v$ in the interior of the region is incident to exactly one boundary
edge, hence we need only show that some $v$ has valency 4.

\bigskip
\leavevmode

\centerline{\epsfbox{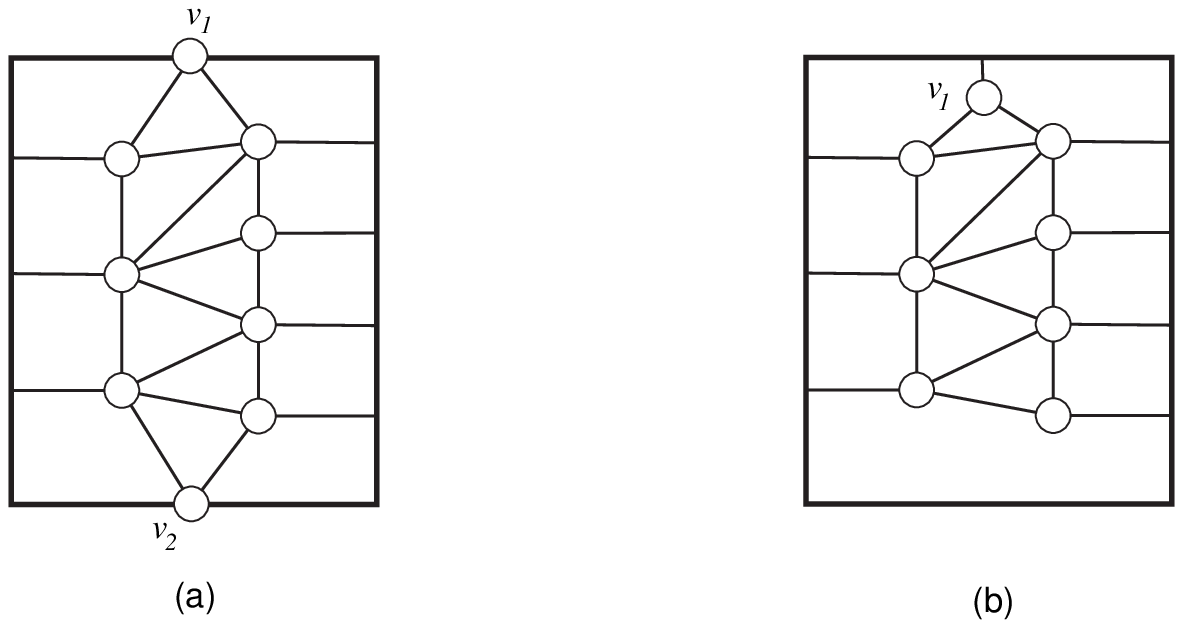}} 
\bigskip
\centerline{Figure 4.2}
\bigskip

If the number of vertices in the interior of the region is 1 or 2, the
result is obvious.  So suppose there are at least 3 such vertices.
Delete $v_2$ and all edges incident to it, and push $v_1$ inwards and
attach a boundary edge to it, as shown in Figure 4.2(b).  Applying
assertion (*) in the proof of Lemma 2.6.5 in [CGLS] to the resulting
graph, we conclude that there is a vertex $v \neq v_1$ of valency at
most 3.  Since there is at most one edge joining $v$ to $v_2$, $v$ has
valency at most 4 (hence exactly 4) in the original graph $\hat G_A$.
\endproof

\proclaim{Proposition 4.8}  Theorem 0.1 is true if $p\geq 3$.
\endproclaim

\proof Let $v$ be the vertex of valency 4 given by Lemma 4.7.  Since
we have assumed $p\geq 3$, by Lemma 2.2 and Proposition 2.4, $G_B$
contains a Scharlemann cycle.  Suppose the Scharlemann cycle has label
pair $\{1, 2\}$.  Then by Lemma 1.2(1), $p$ is even, hence $p\geq 4$,
and the edges of the Scharlemann cycle are not contained in a disk on
$A$.  Thus in $\hat G_A$ there are two edges connecting $v_1$ to
$v_2$, as shown in Figure 5.1(b).  They separate the two boundary
components of the annulus $A$, so no other vertex is incident to two
edges going to different boundary components of $A$.  It follows from
Corollary 4.6 that the only possible vertices of $\hat G_A$ with two
boundary edges are $v_1$ and $v_2$.  Since there are no Scharlemann
cycles on any other label pair, and no extended Scharlemann cycles
(Lemma 1.2(2)), the only labels $x$ for which $\L_x$ has a bigon in
$D$ are $1$ and $2$.  Hence by Lemma 4.2, we have $2p/3 \leq 2$, i.e.\
$p\leq 3$.  But we have just shown that $p\geq 4$, which is a
contradiction.  \endproof

\head \S 5.  The case that $p\leq 2$
\endhead

After Proposition 4.8, it remains to consider the case that the graph
$G_A$ on the annulus $A$ has at most two vertices.  In this section we
will consider this remaining case, and complete the proof of Theorem
0.1.  As before, we assume that $\Delta = 3$.

\proclaim{Lemma 5.1}  If $p\leq 2$, then $p=2$, and the two vertices
of $G_A$ are antiparallel.
\endproclaim

\proof First assume $p=1$.  Then $A$ is a nonseparating annulus in
$M(\a)$.  The reduced graph $\hat G_A$ consists of one vertex, at most
one loop, and at most two boundary edges.  By Lemma 1.1(4) the number
of endpoints of loops is at most $q$.  Since the total valency of
the vertex is $3q$, there exist $q$ parallel boundary edges.  By
Lemma 2.1 $K_{\b}$ is a 1-arch knot.  However, by Proposition 2.4 in
this case $M(\a)$ contains no nonseparating annulus, a contradiction.

Now assume $p=2$ and the two vertices of $G_A$ are parallel.  Then
again $A$ is nonseparating in $M(\a)$.  The reduced graph $\hat G_A$
is a subgraph of one of the two graphs shown in Figure 5.1, depending
on whether or not $\hat G_A$ has a loop.  Since the two vertices are
parallel, by Lemma 1.1(4) each family of parallel interior edges
contains at most $q/2$ edges, hence in both cases there is a family of
at least $q$ parallel boundary edges.  As above, this implies that
$K_{\b}$ is a 1-arch knot, hence contradicts Proposition 2.4 and the
fact that $A$ is nonseparating.  \endproof

\bigskip
\leavevmode

\centerline{\epsfbox{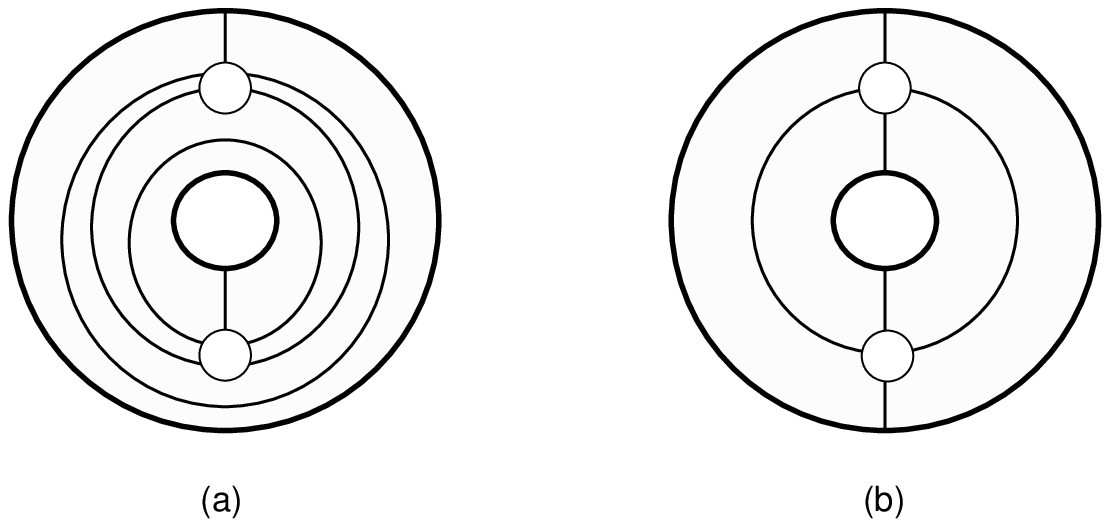}}
\bigskip
\centerline{Figure 5.1}
\bigskip

We may now assume that $p=2$ and the two vertices of $G_A$ are
antiparallel.  Suppose $W$ is a submanifold of $M(\a)$ containing
$K_{\a}$.  We use $\bdd_i W$ to denote the closure of $\bdd W - \bdd
M(\a)$, which is the frontier of $W$ in $M(\a)$, and call it the {\it
interior boundary}.

If $D$ is a disk embedded (improperly) in $M(\a)$ such that $D\cap W =
D \cap \bdd_i W$ is a single arc $C$ on the boundary of $D$, and
$\bdd D - C$ lies on $\bdd M(\a)$, then the pair $(W \cup N(D),
K_{\a})$ is homeomorphic to $(W, K_{\a})$, with $\bdd_i(W\cup N(D))$
identified to $\bdd_i W$ cut along the arc $C$.  This observation will
be useful in the proof of Lemma 5.4.

For the purpose of this section, we define an {\it extremal
component\/} of a subgraph $\Lambda$ of $G_B$ to be a component $\L_0$
such that there is an arc $\gamma$ cutting $B$ into $B_1$ and $B_2$,
with $B_1 \cap \L = \L_0$.

\proclaim{Lemma 5.2}  If $p = 2$ and the two vertices of $G_A$ are
antiparallel, then each vertex of $G_B$ is incident to a boundary
edge.  In particular, each face of $G_B$ is a disk.
\endproclaim

\proof
The reduced graph $\hat G_A$ is a subgraph of that shown in Figure
5.1(a) or (b).  In case (b), each of the interior edge of $\hat G_A$
represents at most $q$ edges of $G_A$, hence each label appear at most
four times at endpoints of interior edges.  It follows that each
vertex of $G_B$ is incident to at least two boundary edges.

In case (a), consider the edge endpoints at a vertex $v$ of 
$G_A$.  Let $s$ be the number of boundary edges at $v$, and let $t$
be the number of loops based at $v$.  Observe that if $s<q$ but
$s+2t > q$ then some label would appear twice among the endpoints of
the parallel loops, which would contradict Lemma 1.1(4).  If $s + 2t
\leq q$, then the two nonloop edges of $\hat G_A$ would represent $3q
- (s+2t) \geq 2q$ edges, which would contradict Lemma 1.1(2).
Therefore we must have $s \geq q$, which implies that each vertex of
$G_B$ is incident to a boundary edge.  

If some face $f$ of $G_B$ is not a disk, then the vertices inside of a
nontrivial loop in $f$ would have no boundary edges, which would
contradict the above conclusion.
\endproof

\proclaim{Lemma 5.3} Suppose $p = 2$ and the two vertices of $G_A$ are
antiparallel.  Then there is a vertex $v_0$ of $\hat G_B$ with the
following properties.

(1)  $v_0$ has valency 2 or 3 in $\hat G_B$, and belongs to a
single boundary edge $e$ of $\hat G_B$.

(2) If the valency of $v_0$ is 3, then the face opposite to the
boundary edge is an interior face.

(3) One of the two faces of $\hat G_B$ containing $e$ intersects $\bdd
B$ in a single arc.
\endproclaim

\proof By Lemma 5.2 each vertex of $G_B$ belongs to a boundary edge.
Consider an extremal component $C$ of $G_B$, and let $\hat C$ be its
reduced graph.  Let $\hat C'$ be its corresponding component in $\hat
G_B$.  Note that $\hat C'$ and $\hat C$ are almost identical, except
that one of the vertices $v'$ of $\hat C'$ may have two parallel
boundary edges, in which case $\hat C$ can be obtained from $\hat C'$
by amalgamating these two edges together.

Note that $\hat C$ must have at least two vertices, for otherwise,
since $C$ is extremal, the vertex would have 6 parallel boundary edges
in $G_B$, two of which would also be parallel on $G_A$ because $\hat
G_A$ has at most four boundary edges, which would then contradict
Lemma 1.1(5).  Hence each vertex of $C$ is incident to at least one
interior edge and one boundary edge, so the valency of each vertex of
$\hat C$ is at least 2.  Modify $\hat C$ as follows.  If some vertex
$v$ of $\hat C$ satisfies condition (1) but not (2), add a boundary
edge to $v$ in the face opposite to the boundary edge at $v$.  Having
done this for all $v$, we get a graph $\hat C''$, which is still a
reduced graph, with at least one boundary edge incident to each
vertex.  Now using (*) in the proof of [CGLS, Lemma 2.6.5] and arguing
directly when $\hat C''$ has only two or three vertices, we see that
$\hat C''$ contains at least two vertices, each of which has valency 2
or 3 in $\hat C''$ and belongs to a single boundary edge of $\hat
C''$.  At least one of these two vertices, say $v_0$, is not the
vertex $v'$ above, hence it has property (1) when considered as a
vertex in $\hat G_B$.  By the definition of $\hat C''$, $v_0$
automatically has property (2).  To prove (3), notice that if both
faces containing $e$ intersect $\bdd B$ in more than one arc, then $C$
would not be an extremal component.  \endproof

\proclaim{Lemma 5.4} If $p = 2$, and the two vertices of $G_A$ are
antiparallel, then $\bdd M$ is a union of tori.  \endproclaim

\proof Let $W_0$ be a regular neighborhood of $A \cup K_{\a}$.  Since
$K_{\a}$ intersects $A$ in two points of different signs, $\bdd_i W_0$
has two components $F_b, F_w$, each being a twice punctured torus.  The
annulus $A$ cuts $W_0$ into two components $W_0^b$ and $W_0^w$ (with
$W_0^b \supset F_b$), which will be called the black region and the
white region, respectively.  If $E$ is a disk face of $G_B$ or more
generally a disk in $M(\a)$, then $E$ is said to be {\it black\/}
(resp.\ {\it white}) if $E \cap W_0$ lies in the black (resp.\ white)
region.

Suppose $D$ is a compressing disk of $F_b$ in $M(\a) - \Int W_0$.  If
$\bdd D$ is a nonseparating curve on $F_b$, then after adding the
2-handle $N(D)$ to $W_0$, the surface $F_b$ becomes an annulus.  If
$\bdd D$ is separating on $F_b$, then it is not parallel to a boundary
curve of $F_b$ because $\bdd F_b$ is parallel to $\bdd A$ and $A$ is
incompressible in $M(\a)$; thus $\bdd D$ must cut $F_b$ into a once
punctured torus and a thrice punctured sphere, and after adding the
2-handle $N(D)$ the surface $F_b$ becomes the union of a torus $S_1$
and an annulus $S_2$.  Since $M$ is simple, $S_1$ either is boundary
parallel or bounds a solid torus, and $S_2$ must be boundary parallel,
because $\bdd S_2$ is parallel to $\bdd A$ and $A$ is incompressible,
which implies that $S_2$ is incompressible.  In any case, we have
shown that if $F_b$ is compressible in $M(\a) - \Int W_0$ then there
is a component $C_b$ of $M(\a) - \Int W_0$ such that $C_b \cap W_0 =
F_b$ and $C_b \cap \bdd M(\a)$ is either an annulus or the union of an
annulus and a torus.  Similarly for $F_w$.  In particular, if both
$F_b$ and $F_w$ are compressible in $M(\a) - \Int W_0$, then $\bdd
M(\a)$ is a union of tori, and we are done.  From now on, we will
assume that $F_w$ is incompressible in $M(\a) - \Int W_0$ and show
that this will lead to a contradiction.  Note that the assumption
implies that $G_B$ has no interior white face: For, by Lemma 5.2 all
faces of $G_B$ are disks, and since $G_A$ has no trivial loops, the
boundary of an interior face is always an essential curve on $\bdd_i
W_0$; hence an interior white face would give rise to a compressing
disk of $F_w$ in $M(\a) - \Int W_0$.

Let $v_0$ be a vertex of $G_B$ given by Lemma 5.3.  By Lemma 5.3(1),
$v_0$ has valency 2 or 3 in $\hat G_B$.  First assume that $v_0$ has
valency 2 in $\hat G_B$.  Then the interior edge $e$ of $\hat
G_B$ incident to $v_0$ must represent exactly two edges of $G_B$: It
cannot represent more than two edges, otherwise there would be two
interior faces of different colors, contradicting the fact that
$G_B$ has no interior white face.  It cannot represent only one edge
of $G_B$, otherwise $v_0$ would have five parallel boundary edges,
which would contradict Lemma 1.1(5) because $\hat G_A$ has at most
four boundary edges.  Thus the part of $G_B$ near $v_0$ is as shown in
Figure 5.2, where $f$ is the interior (black) face bounded by the two
edges represented by $e$.

Now assume that $v_0$ has valency 3 in $\hat G_B$.  Then by Lemma
5.3(2) the face $f$ of $\hat G_B$ opposite to the boundary edge at
$v_0$ is an interior face.  Thus $f$ is a black face, and each of the
interior edges of $\hat G_B$ incident to $v_0$ represents only one edge
of $G_B$ as otherwise there would be a white interior face.  Hence
again the part of $G_B$ near $v_0$ is as shown in Figure 5.2.

Consider the white boundary faces $D_0, D_1, D_2$ as shown in Figure
5.2.  By Lemma 5.3(3), we may assume that $D_1$ intersects $\bdd B$ in
a single arc.  Let $C_i$, $i = 0, 1$, be the arc $D_i \cap F_w$.  Then
$C_0, C_1$ are essential arcs on $F_w$.  Moreover, since $C_0$
intersects a meridian of $K_{\a}$ exactly once, while $C_1$ intersects
it at least twice, they are nonparallel.  Recall that $\bdd _i (W_0
\cup N(D_0 \cup D_1))$ is obtained from $\bdd _i W_0$ by cutting along
$C_0\cup C_1$.  Since $C_0$ and $C_1$ are nonparallel, they cut $F_w$
into one or two annuli, which must be boundary parallel because we
have assumed that $F_w$ is incompressible and $M$ is simple.  It
follows that the whole surface $F_w$ is boundary parallel.  Now a
meridian disk of $N(K_{\a})$ in the white region corresponds to a disk
in $M(\a)$ intersecting the curve $K_{\a}$ in a single point, which
gives rise to an essential annulus in $M$, contradicting the fact that
$M$ is anannular.  \endproof

\bigskip
\leavevmode

\centerline{\epsfbox{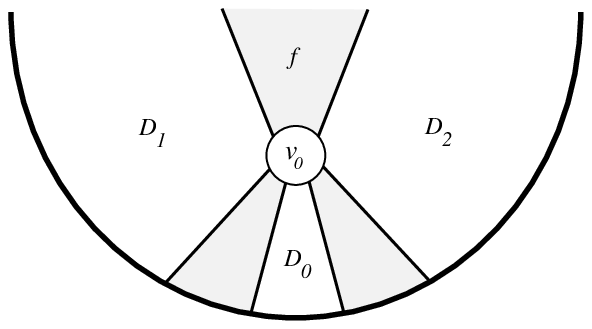}}
\bigskip
\centerline{Figure 5.2}
\bigskip

\noindent
{\it Proof of Theorem 0.1.}  By Proposition 4.8, we may assume that
$p\leq 2$.  By Lemmas 5.1 and 5.4, $\bdd M$ is a union of tori.  Since
$M(\b)$ is $\bdd$-reducible, either it is reducible or it is a solid
torus.  In the first case the result follows from [Wu3, Theorem 5.1].
So we assume that $M(\b)$ is a solid torus.  In particular, $\bdd
M(\a)$ is a single torus $T$.  The boundary of the annulus $A$ cuts
$T$ into two annuli $A_1, A_2$.  If some $A\cup A_i$ is an essential
torus in $M(\a)$ then $M(\a)$ is toroidal, so the result follows from
[GL4].  If each $A\cup A_i$ is inessential, then it bounds a solid
torus (note that it cannot be boundary parallel, otherwise $A$ would
be boundary parallel).  It follows that $M(\a)$ is a Seifert fiber
space with orbifold a disk with two singular points.  It was shown in
[MM1, Theorem 1.2] that if $M(\a)$ is a Seifert fiber space and $M(\b)$
is a solid torus then $\Delta \leq 1$.  This completes the proof of
Theorem 0.1.  
\qed

In the proof of Theorem 0.1, we assumed that the manifold $M$ is
simple.  However, the conditions that $M$ is irreducible and atoroidal
can be removed from the assumptions.

\proclaim{Corollary 5.5} Suppose $M$ is anannular and boundary
irreducible.  If $M(\a)$ is annular and $M(\b)$ is boundary reducible,
then $\Delta(\a, \b) \leq 2$.  \endproclaim

\proof First assume that $M$ is irreducible but toroidal.  Since $M$
is anannular, by the canonical splitting theorem of
Jaco-Shalen-Johannson (see [JS, p.\ 157]) there is a set of essential
tori $\Cal T$ cutting $M$ into a manifold $M'$ such that each
component of $M'$ is either a Seifert fiber space or a simple
manifold.  If the component $X$ containing the boundary torus
$\partial_0 M$ is Seifert fibered, then it contains an essential
annulus consisting of Seifert fibers, with both boundary components on
$\partial_0 M$, so $M$ would be annular, contradicting our assumption.
So assume $X$ is simple.  Since $M(\b)$ is boundary reducible, by
looking at a boundary reducing disk $B$ which has minimal intersection
with $\Cal T$, one can see that $X(\b)$ must be boundary reducible.
Similarly one can show that $X(\a)$ is either boundary reducible or
annular.  Applying Theorem 0.1 and [Wu2, Theorem 1] to $X$, we have
$\Delta \leq 2$.

If $M$ is reducible, split along a maximal set of reducing spheres to
get an irreducible manifold $M'$.  By an innermost circle argument one
can show that $M'(\a)$ is annular and $M'(\b)$ is boundary reducible,
so the result follows from that for irreducible manifolds.
\endproof

\enddocument